# The complete vertex *p*-center problem

F. Antonio Medrano[1]

**Abstract.** The vertex *p*-center problem consists of locating *p* facilities among a set of *M* potential sites such that the maximum distance from any demand to its closest located facility is minimized. The complete vertex *p*-center problem solves the *p*-center problem for all *p* from 1 to the total number of sites, resulting in a multi-objective trade-off curve between the number of facilities and the service distance required to achieve full coverage. This trade-off provides a reference to planners and decision-makers, enabling them to easily visualize the consequences of choosing different coverage design criteria for the given spatial configuration of the problem. We present two fast algorithms for solving the complete *p*-center problem, one using the classical formulation but trimming variables while still maintaining optimality, the other converting the problem to a location set covering problem and solving for all distances in the distance matrix. We also discuss scenarios where it makes sense to solve the problem via brute-force enumeration. All methods result in significant speed-ups, with the set covering method reducing computation times by many orders of magnitude.

## 1 Introduction

Let *N* be the set of discrete demands to be served, represented by the index $i = 1, 2, \ldots, n$. Let *M* be the set of discrete potential facility sites that may serve those demands, represented by the index $j = 1, 2, \ldots, m$. Let $d_{ij}$ be the distance from demand *i* to site *j*. The vertex *p*-center problem consists of locating *p* facilities among a set of *m* potential

[1]  F. Antonio Medrano
    antonio.medrano@tamucc.edu

    Conrad Blucher Institute for Surveying and Science, Department of Computing Sciences, Texas A&M University–Corpus Christi, 6300 Ocean Drive, Unit 5799, Corpus Christi, TX, USA

sites such that the maximum distance from any demand to its closest located facility is minimized. This article considers the problem in the Euclidean space, where the "service distance" is a real number radius of coverage around the located facilities.

The *p*-center (PC) problem was originally formulated by (Hakimi 1964) and it is often used for equitable location of public services such as fire stations, police stations, ambulance depots (Daskin 2013; Marianov and ReVelle 1995; Shier 1977) since it minimizes the furthest distance of the service to any constituent. Variations include locating facilities anywhere in the Euclidean space (Drezner 1984), also known as the continuous or non-vertex problem; locating facilities on a network (Kariv and Hakimi 1979); covering continuous area demands rather than point demands (Suzuki and Drezner 1996); and partial coverage variants (Daskin and Owen 1999). Many other variations exist, as thoroughly cataloged in Calik et al. (2015) and Contardo et al. (2019). The formulations in this article expand upon the vertex *p*-center problem with discrete point facility site locations and demands, and it can be assumed henceforth that we always refer to this vertex *p*-center problem as opposed to any other variant.

The *p*-center problem is strongly NP-Hard, since its formulation is of the form of a location-allocation mixed-integer programming problem. Location-allocation problems are typically not integer friendly (ReVelle 1993), and often require relaxations or heuristics to solve for large problem sizes (Dzator and Dzator 2013; Gadegaard et al. 2018; Mladenović et al. 2003; Rosing et al. 1979; Teitz and Bart 1968). Minieka (1970) was the first to realize the close relationship between the *p*-center problem and the location set covering problem (LSCP), where the LSCP consists of minimizing the number of facilities to completely cover a set of demands given each facility has a fixed covering radius. He used this relationship to heuristically solve the *p*-center problem by solving numerous LSCP problems, adjusting the coverage distance until a suitable bound was found. More recently, Elloumi et al. (2004) used the fact that optimal solutions to the *p*-center problem will always be distances that exist in the distance matrix for the problem (Daskin 2013). They developed a clever two-stage formulation that first converts the problem from a MIP to a binary radius formulation problem by replacing the continuous covering radius variable with a set of binary variables that represent "coverage disks" with radii of unique distances $D^k$ that exist in the distance matrix. They then used an LP relaxation of the formulation to determine distance lower and upper bounds for a given *p*, and iterated an LSCP in-between the bounds to optimality. This formulation is very complex, but harnesses numerous deep insights into the *p*-center problem to create a faster solution algorithm.

Further refinements to finding tighter bounds were proposed by Calik and Tansel (2013), who presented a method similar to that of Elloumi et al. (2004) but at each iteration could update both the upper and lower bounds, thus converging on the solution more quickly. The most recent work on solving the PC problem comes from Contardo et al. (2019), who developed an improvement on the relaxation-based iterative algorithm by Chen and Chen (2009). The relaxation approach considers only local subsets of the demands within the coverage matrix. It assigns demands to facilities in the local subsets, and if all demands are allocated then the solution is globally optimal. If some demands are not allocated to any facility, then the local subsets are expanded with a fixed number of demands, and the process repeats until all demands are assigned and the optimal solution is found. The method by Contardo et al. (2019), which they call CIK2018, is memory efficient by avoiding the computation and storage of the complete covering matrix, allowing it to solve for single *p*-center problems of up to one million points, and more quickly than previous methods. Their algorithm was developed for integer distance matrices and for small values of *p*, from 2 to 30 in their paper.

Some studies have modeled bi-objective center problems where the number of facilities *p* and the covering distance are varied. One such model, by Niblett and Church (2015), is intended for retail site location. It trades-off the number of facilities with the separation distance between facilities in order to maximize retail profits, while adding some amount of over-coverage in order to prevent competition from other franchises to profitably enter a market. Other bi-objective franchise location models that use *p*-center for one objective include a greedy metaheuristic by Salazar-Aguilar et al. (2013), and the *p*-center *p*-dispersion bi-criteria model by Tutunchi and Fathi (2019).

## 1.1 Relationship between covering problems

Daskin (2013) outlines the close relationship that exists between the three primary types of covering location problems, namely the PC problem, the LSCP, and the maximal covering location problem (MCLP) (Church and ReVelle 1974). All three of these contain elements of the amount of coverage, the number of facilities, and the facility covering distance; in which one sets two and optimizes for the third. The relationship is as follows.

1. **The *p*-Center Problem**
   - Fixed number of facilities and achieve complete coverage.
   - Minimize the coverage distance.

2. **The Location Set Covering Problem**
   - Fixed coverage distance and achieve complete coverage.
   - Minimize the number of facilities.

3. **The Maximal Covering Location Problem**
   - Fixed number of facilities and coverage distance.
   - Maximize the amount of coverage.

The MCLP typically is solved in instances when the number of facilities and the coverage distance cannot achieve complete coverage. As such, it maximizes the weighted demand that can be covered. The close relationship between these three covering problems is clear when plotting the results of the complete *p*-center problem.

## 2 Problem Formulation

The complete *p*-center problem solves the *p*-center problem for all *p* from 1 to *m*. Mathematically, this can be formulated with the classical mixed-integer programming formulation for the *p*-center problem to create a model we call CPC. Let $y_j = 1$ if a facility is located at site *j*, 0 if not. Let $x_{ij} = 1$ if demand *i* is allocated to facility *j*, 0 if not. Let $z_p$ be the minimum coverage radius require to cover all demands with *p* facilities. We already defined $d_{ij}$ to be the distance from demand *i* to site *j*. Then the CPC is the following.

**Formulation (CPC)**

$$\textbf{for } p = 1, \ldots, m;\ \text{minimize } z_p \tag{1}$$

subject to

$$\sum_{j=1}^{m} y_j \leq p \tag{2}$$

$$\sum_{j=1}^{m} x_{ij} = 1 \qquad i = 1, \ldots, n \tag{3}$$

$$x_{ij} \leq y_j \qquad i = 1, \ldots, n;\ j = 1, \ldots, m \tag{4}$$

$$\sum_{j=1}^{m} d_{ij} x_{ij} \leq z_p \qquad i = 1, \ldots, n \tag{5}$$

$$x_{ij}, y_j \in \{0,1\} \qquad i = 1, \ldots, n;\ j = 1, \ldots, m \tag{6}$$

The for–loop in constraint (1) is the only difference between the CPC and the classical $p$-center formulation, where in each iteration it sets the value of $p$ in constraint (2) to be the number of facilities to be located. Constraint (3) ensures that each demand is allocated to one facility. Constraint (4) is a set of Balinski constraints that ensure demands are only allocated to located facilities (Balinski 1965). Constraints (2), (3) and (4) make this formulation a location-allocation type problem (Church and Medrano 2018). Constraint (5) defines the value of the coverage distance that is minimized in the objective, and (6) defines the binary variables. This formulation, for each iteration of $p$ in the for–loop, has $nm + m$ binary variables, 1 continuous variable, and $nm + 2n + 1$ constraints. In our computational experiments, problems solved with this formulation are denoted CPC-MIP.

When solving the CPC problem, the solution is the set $Z = \{z_1, z_2, \ldots, z_m\}$, corresponding to the number of facilities located $p = 1, 2, \ldots, m$. For simplicity, all experiments in this article solve the symmetric variant where every demand location is also a potential facility site. The resulting solution is a trade-off curve between the number of facilities and the minimum covering distance required to achieve complete coverage, as shown in Fig. 1 for the 55-point Swain data set (Swain 1971) commonly used in location optimization literature (Church and Weaver 1986; Church and Gerrard 2003; Church and Roberts 1983; Daskin and Owen 1999; Dzator and Dzator 2013; Schilling 1982; Snyder and Daskin 2006; Solanki 1991). Above the CPC curve is the set covering region, where one may find numerous configurations to achieve complete coverage. Below the CPC curve is maximum coverage region, where one cannot achieve complete coverage with the given $p$ and coverage radius, and must solve an MCLP in order to maximize the coverage.

With this figure in hand, one can instantly know by choosing any particular number of facilities and coverage radius whether one can completely cover all demands with those parameters. Or the designer can evaluate how many facilities will be required to achieve complete coverage at different service distance standards. Or if the designer has different budget options for the number of facilities they can build, they can instantly evaluate what the service distance would be for those budgets. This is a simple yet powerful design tool for decision makers to use when debating alternative plans, that provides solutions to all three primary types of covering problems with one graphic. Fo and da Silva Mota (2012) ran a case study to locate healthcare facilities in Brazil, where they compared results from running a *p*-median model, a *p*-center model, a set covering model, and a maximal covering model. Having a CPC curve such as the one in Fig. 1 would have provided information on three out of their four models. If one were to solve the partial cover *p*-center problem (Daskin and Owen 1999) completely, one could potentially incorporate isolines of partial coverage to this figure. For now, we leave this possibility open to future research.

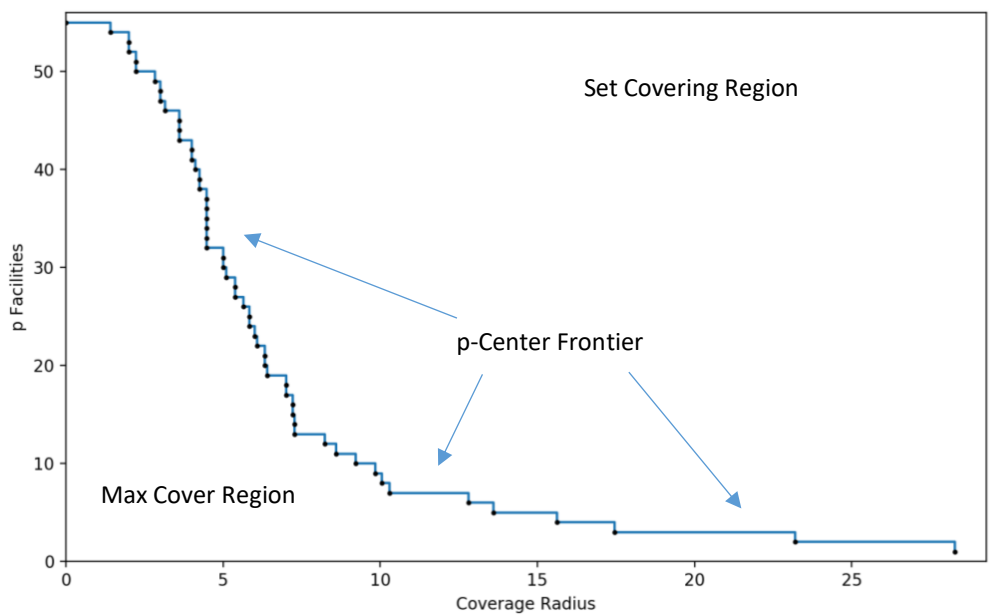

**Fig. 1.** The CPC solution to the 55-point Swain data set, displaying the frontier of the minimum covering radius required with different numbers of facilities.

For OR scientists, the CPC problem also provides motivation to develop *p*-center methods that are efficient for all *p*-values. In order to perform well, a CPC problem algorithm must be fast from $p = 1, 2, \ldots, m$. The current state-of-the-art for the integer PC problem, CIK2018 (Contardo et al. 2019), is tested for up to $p = 30$ on data sets with up to 1 million points. It is easy to foresee applications on big data where there is a need to locate a greater number of facilities, thus an efficient CPC problem method will be able to quickly solve for all *p*-values of interest. We provide ideas at the end of this paper for additional methods that use a combination of algorithms to solve the CPC problem efficiently for all *p*-values on big data.

## 3 Solution Methodologies

### 3.1 Trivial Solutions

The CPC problem has a few trivial solutions, namely when $p = 1$, when $p = m - 1$, and when $p = m$. Let $D^k$ be a sorted list of the $k = 0, …, K$ unique distances that exist in the symmetric distance matrix, where $D^0 = 0$ and $D^K$ is the maximum distance between any two points in the data set. When $p = m$, the solution is $z_m = 0$, since every facility can cover itself. When $p = m - 1$, the solution is $z_{m-1} = D^1$, where the two closest demands are covered by one facility with a covering radius that is the smallest distance in the distance set. When $p = 1$, the solution $z_1$ is the minimum row element of the maximum column distance in the distance matrix, essentially the facility with the smallest maximum distance to any demand. These operations are essentially instantaneous to compute for any sized problem, and should always be used for such $p$.

### 3.2 Real Number versus Integer Distances

This paper solves the CPC problem for Euclidean real number distances. Much of the literature focuses instead on using integer distances, with the caveat that decimal point precision could be incorporated via scaling. It is important to note that rounding distances to integers adds a quantization error to the distance matrix, affecting the number of unique distances in the matrix based on the range of the data. If the data has Euclidean distances that range from [1,10], then there will be a potential of 10 unique distances in the distance matrix after rounding. If the same data is then scaled to [10,100] before rounding, then there are 91 potential unique distances in the distance matrix, and so on. Every application should specify its tolerance for distance precision to ensure accurate problem representation.

In this paper, we solved all problems for Euclidean real number distances. In section 4.3 we compare our methods on the same problems with both integer and real number distances, and show how the real number variant is overall much harder to solve. Using scaling to solve to some number of decimal points will result in a trade-off between problem difficulty and precision. In this article, algorithms solving the complete $p$-center problem for real numbers are denoted with the CPC prefix, and algorithms for the integer variant are denoted with CPCi. Our integer distance trials round all distances up, ensuring that all solutions to the CPCi are feasible integral upper bounds to their associated CPC real number solutions for all $p$-values.

### 3.3 MIP with Upper Bound Variable Reduction (CPC-MIPUB)

Variable reduction of a location-allocation mathematical optimization model while still maintaining guaranteed optimality has been shown to be a very effective way to speed up its computation (Church 2008; Church 2016; García et al. 2011). The key to the complete $p$-center mixed-integer programming with upper bound model (CPC-MIPUB) is to find variables associated with distances in the distance matrix that cannot represent an optimal solution, and eliminate them from the problem. Since the trade-off curve is monotonically decreasing with respect to $z_p$ when iterating to larger $p$-values, this means that for any $p + 1$, the objective $z_{p+1} \leq z_p$. Since $z_p$ decreases as $p$ increases, any variables associated with distances greater than $z_p$ can be eliminated from the problem for any future larger *p* while still maintaining optimality. This can be mathematically formulated as such, letting $UB = z_{p-1}$ for any *p*-value.

**Formulation (CPC-MIPUB)**

$$\textbf{for } p = 1, \ldots, m; \ \text{minimize } z_p \tag{7}$$

subject to

$$\sum_{j=1}^{m} y_j \leq p \tag{8}$$

$$\sum_{j:d_{ij} \leq UB}^{m} x_{ij} = 1 \quad i = 1, \ldots, n \tag{9}$$

$$x_{ij} \leq y_j \quad i = 1, \ldots, n; \ j: d_{ij} \leq UB \tag{10}$$

$$\sum_{j:d_{ij} \leq UB}^{m} d_{ij} x_{ij} \leq z_p \quad i = 1, \ldots, n \tag{11}$$

$$x_{ij}, y_j \in \{0,1\} \quad i = 1, \ldots, n; \ j: d_{ij} \leq UB \tag{12}$$

Constraints (7) and (8) are the same as for the CPC-MIP. Constraint set (10) has fewer Balinsky constraints, since one eliminates any for allocations of a distance greater than the *UB*. Constraints (9), (11), and (12) eliminate any $x_{ij}$ variables associated with distances greater than *UB*.

Given that finding $z_1$ is trivially fast, one could always speed up the computation for a single *p*-value using the MIP model by first solving for $p = 1$ and then trimming variables for distances larger than $z_1$. After presenting the CMP-MIPUB formulation at a conference, the author became aware that the same reduction technique had been developed for the single *p*-center problem in a dissertation by Wei (2008) but never published in a journal. They iterated a few guesses for an upper bound distance, rather than using a previously found solution as the cutoff.

### 3.4 Iterative LSCP Algorithm (CPC-LCSP)

Minieka (1970) discovered the close relationship between the PC problem and the LSCP. Elloumi et al. (2004) used this close relationship, and the knowledge that every solution to the PC problem will be a unique distance that exists in the distance matrix, to quickly solve for a single fixed *p* solution. Recall the LSCP uses a fixed covering distance to achieve complete coverage with the minimal number of facilities. We propose a method for solving the CPC where we iterate the LSCP for all unique distances in the distance matrix until the complete solution set has been found. While the MIP formulation solves the CPC in *m* iterations, the LSCP solves in O(*nm*) iterations, but each iteration of the LSCP problem solves much more quickly than the MIP formulation. The MIP formulation has $nm + m$ binary variables, 1 continuous variable, and $nm + 2n + 1$ constraints, and is not integer friendly. The LSCP has *m* variables and $n$ constraints, and is integer friendly. Being integer friendly means that solving the LP relaxation will often give an integer solution, thus the problem can be solved much more quickly by first trying an LP solver and seeing if the solution is integral, rather than using branch-and-bound and cutting plane methods used in IP solvers. While there are many more iterations with CPC-LCSP than with CPC-MIP, each iteration solves much more quickly.

The CPC-LSCP method starts with $z_m = D^0 = 0$ and $z_{m-1} = D^1$, since these solutions are trivial. It then increments over $k$ in the set $D^k$, setting each value as the new covering distance, and solves the LSCP to determine which distances trigger a decrease in the number of facilities needed for complete coverage, thus solving for $z_{m-2}, z_{m-3}, \ldots, z_2$. Finally, one uses simple matrix operations to find $z_1$ in order to complete the solution set. For each iteration, one must update the coverage matrix $A$ with the new covering distance $D^k$, then update the model constraints accordingly.

**Algorithm (CPC-LSCP)**

    Solve for $z_1$, $z_m$ and $z_{m-1}$ trivially.
    $p = m–2, \; k = 2$

    **while** ( $p \geq 2$ ) {
        Update the coverage matrix $A$ for the covering distance $D^k$.
        Solve (LSCP) with the current $A$ to compute $z_{LSCP}$.

        **while** ( $z_{LSCP} < p$ ) {
            $z_p = D^k$
            $p = p - 1$
        }
        $k = k + 1$
    }

The generic formulation for the location set covering problem is as follows.

**Formulation (LSCP)**

$$\text{minimize } z_{LSCP} = \sum_{j=1}^{m} y_j \tag{13}$$

subject to

$$\sum_{j=1}^{m} a_{ij} y_j \geq 1 \qquad \text{for } i = 1, \ldots, n \tag{14}$$

$$y_j \in \{0,1\} \qquad \text{for } j = 1, \ldots, m \tag{15}$$

CPC-LSCP is particularly fast for high $p$-values, since the covering distance gap between successive values tends to be smaller for high values than for lower values. CPC-LSCP also take advantage of the property that numerous high $p$-values may have the same optimal service distance solution. I will present an example here, and Daskin (2013) elaborates on this in the first section of his book chapter on the *p*-center problem. In Fig. 1, one can see for the Swain data trade-off graph a vertical sequence of points where $z_{38} = 4.242641$, and $z_{37} = z_{36} = z_{35} = z_{34} = z_{33} = z_{32} = 4.472136$. This was computed by first evaluating LSCP with covering distance 4.242641 and getting a solution of 38, then iterating to the next unique covering distance, 4.472136, and getting a solution of 32. Essentially, two solver iterations find seven solutions because distance ties in the distance matrix cover many demands with the next larger unique distance. The Swain data set is very small, and for larger data sets CPC-LSCP can solve for large sections of the high *p*-range with very few iterations.

We tested using the row, column, and essential site domination rules in Church and Gerrard (2003), which use linear algebra to reduce the covering matrix problem size while maintaining optimality. While these did successfully speed up the code on slower open-source solvers we tested, we found that that with Gurobi the domination computations actually slowed down finding solutions for the problem sizes that we evaluated (up to 1,400 points). This is due to the sophisticated pre-solvers Gurobi has, which likely implement many of these domination rules using highly-optimized code. For much larger problem sizes, where memory becomes an issue, implementation of the domination rules could provide a performance speedup.

### 3.5 Iterative LSCP Algorithm with Brute-Force Enumeration (CPC-LSCPe)

All combinatorial optimization problems could be solved by enumerating every possible solution, evaluating how each one performs, and choosing the best performing as the optimal solution. Typically, the combinatorics explode factorially with problem size making this approach unreasonable for all but the smallest problems. For the CPC, there are a very few cases where enumeration may yield faster results.

We already used simple linear algebra to compute $z_1$, but one could instead try evaluating each site, seeing what the maximum demand distance is in the distance matrix for locating a facility at that site, and picking the site with the smallest maximum distance as the solution. This would require *m* iterations of picking and checking.

For two facilities, this is simple as well. If one picks sites 1 and 2, the columns in the distance matrix corresponding to these are $d_{i1}$ and $d_{i2}$. The minimum distance to achieve complete coverage is the maximum of the element-wise minimum of these two columns, $\max(\min(d_{i1}, d_{i2}))$. With $\frac{1}{2}(m^2 - m)$ iterations in a nested for-loop, one could check every possible combination to find $z_2$ for the PC problem. Similarly, for three sites one takes the maximum of the element-wise minimum of the three columns in three nested for-loops, and so forth. In general, the combinatoric iterations follow the rule

$$\binom{m}{p} = \frac{m!}{(m-p)!\,p!} \qquad (16)$$

This has an exponential complexity of $O(m^p)$, which is unreasonable for large *p*-values. But for small *p*, if one can make the enumeration code run extremely fast, enumeration is the simplest and potentially the fastest way to compute $z_p$. Our experimental code implements everything in Python, which is by no means the fastest language for repeated basic math calculations. But for enumeration we used the Numba package, which takes highly numerical Python code snippets and uses LLVM to compile them into much faster machine code. Numba can approach speeds of C or Fortran, and has parallel for-loops to take full advantage of multi-core processors. PCP-LSCPe using enumeration for $p = 2$ and 3 always yielded faster results with real number distances than without the enumeration. Otherwise for these two *p-values*, PCP-LSCP has to iterate over many unique distances in the distance matrix to determine which is the cutoff. For p > 3, computation times with enumeration were slower than for using enumeration with only $p = 2$ and 3, although some hardware such as 20+ core workstations may be capable of enumeration speedups with $p = 4$. For our 6-core hardware enumeration up to $p = 3$ was the sweet spot for speed gains, so the CPC-LSCPe methods in this article use enumeration only for $p = 2$ and 3 and use CPC-LSCP for all other *p*-values. The Github repository with our code does contain a CPC-LSCPe implementation for $p \leq 4$, and test files for $p = 2, 3, 4$, and 5.

## 4 Computational Experiments

We tested the efficiency of the various CPC algorithms on a variety of spatial data sets. For simplicity, the experiments here solve the symmetric variant where every demand location is also a potential facility site, i.e. $N = M$. All derived distances are Euclidean real number distances unless noted otherwise.

### 4.1 Implementation Details

All models were tested with Gurobi 9.0.1 using their Python 3.7.7 interface. Enumeration was done using the Numba 0.42.1 just-in-time LLVM compiler, which compiles and parallelizes sections of the Python code. All computations were run on an Apple MacBook Pro computer with a six-core Intel Core i9-8950HK processor, 32GB of RAM, and running macOS 10.14. The reader is invited to examine and download all codes from Github (Medrano 2019).

The Github repository also contains the data sets used for testing. It contains the 55-point Swain data set (Swain 1971), a tiny data set used primarily for prototyping and development. It also contains the TSPLIB Traveling Salesman Problem Library (Reinelt 1991) problems, used for more comprehensive benchmarking.

Let us mention a couple other minor implementation notes on our code that is available for download. First, or any *p*-value the CPC-MIP, CPC-LSCP, and CPC-LSCPe formulations have the same number of variables and constraints, thus the Gurobi solver implicitly applies warm starts to each iteration to speed up computation. For CPC-MIPUB, the model changes enough by deleting variables and constraints at each iteration that the solver does not automatically use warm starts. Instead, our code explicitly applies warm starts to CPC-MIPUB, which speeds up computations by about 27% compared to without them, while still maintaining optimality. Second, for the CPC with real numbers, rather than a Euclidean distance matrix, one can instead run the problem with a squared-Euclidean distance matrix and then take the square root of each solution for each *p*-value. This saves a bit of time in calculating the distance matrix since square root operations are computationally slow, saving $O(n^2)$ square root operations in computing the distances, then computing $O(n)$ square roots to output the solution. Doing so is only valid because the distances are used to determine single-level coverage, distances are not used as part of a cost-distance function as in the *p*-median problem (Daskin and Maass 2015). We have confirmed that for real number distances that this squared-Euclidean method provides identical results to computing with the Euclidean distance matrix, while reducing the time spent on calculating the distance matrix. For integer distance matrices though this is not a valid approach, since rounding on a squared distance results in non-linear quantization. Some integer distance methods for computing complete $p$-center problem solutions are faster for other reasons, but real number approaches can save a bit of time with this trick.

**Table 1** Aggregated results of computational experiments comparing the four algorithms for solving the CPC problem with real number distances.

| Name | Size | CPC-MIP Time (s) | CPC-MIPUB Time (s) | CPC-MIPUB Speedup | CPC-LSCP Iter. | CPC-LSCP Time (s) | CPC-LSCP Speedup | CPC-LSCPe Iter. | CPC-LSCPe Time (s) | CPC-LSCPe Speedup |
|---|---|---|---|---|---|---|---|---|---|---|
| swain.dat | 55 | 2.85 | 1.64 | 1.74 | 183 | 0.15 | 18.43 | 93 | 0.11 | 27.01 |
| st70.tsp | 70 | 5.81 | 2.79 | 2.08 | 572 | 0.54 | 10.68 | 270 | 0.21 | 28.05 |
| rd100.tsp | 100 | 23.04 | 9.72 | 2.37 | 2347 | 4.00 | 5.76 | 1246 | 1.82 | 12.66 |
| bier127.tsp | 127 | 59.00 | 25.94 | 2.27 | 1417 | 3.65 | 16.16 | 692 | 2.17 | 27.22 |
| u159.tsp | 159 | 95.34 | 26.40 | 3.61 | 235 | 0.81 | 117.88 | 105 | 0.74 | 129.14 |
| rat195.tsp | 195 | 367.12 | 190.97 | 1.92 | 2190 | 16.80 | 21.85 | 1330 | 11.71 | 31.35 |
| d198.tsp | 198 | 237.93 | 74.84 | 3.18 | 4867 | 18.89 | 12.60 | 2120 | 7.61 | 31.28 |
| gr202.tsp | 202 | 331.33 | 160.70 | 2.06 | 19114 | 127.09 | 2.61 | 15940 | 108.31 | 3.06 |
| tsp225.tsp | 225 | 449.54 | 184.13 | 2.44 | 5907 | 48.68 | 9.23 | 3451 | 25.96 | 17.32 |
| gil262.tsp | 262 | 853.57 | 355.52 | 2.40 | 3067 | 38.27 | 22.30 | 1260 | 15.62 | 54.66 |
| a280.tsp | 280 | 1196.91 | 558.63 | 2.14 | 612 | 8.27 | 144.76 | 351 | 6.36 | 188.12 |
| pr299.tsp | 299 | 1414.26 | 604.88 | 2.34 | 4670 | 62.38 | 22.67 | 2551 | 34.85 | 40.58 |
| lin318.tsp | 318 | 1653.83 | 443.21 | 3.73 | 7698 | 113.64 | 14.55 | 4032 | 50.40 | 32.81 |
| gr431.tsp | 431 | 6669.41 | 2727.10 | 2.45 | 73853 | 1491.93 | 4.47 | 49970 | 973.80 | 6.85 |
| pr439.tsp | 439 | 6687.94 | 1560.88 | 4.28 | 10071 | 234.47 | 28.52 | 4782 | 103.72 | 64.48 |
| u574.tsp | 574 | 27877.83 | 10910.72 | 2.56 | 64112 | 3694.59 | 7.55 | 45098 | 2453.82 | 11.36 |
| p654.tsp | 654 | 39214.96 | 3838.07 | 10.22 | 5158 | 106.02 | 369.89 | 2952 | 58.44 | 671.05 |

## 4.2 Experimental Evaluation with Real Number Distance

The CPC-MIP, CPC-MIPUB, CPC-LSCP, and CPC-LSCPe implementations were all tested on the Swain and the TSPLIB spatial data sets using real number distances. Table 1 shows the experimental computation results. First it shows the dataset name and the number of points in the dataset. All given speedups are relative to the CPC-MIP times. Iterations are not listed for the MIP implementations because they have the same number of iterations as the size of the dataset. For CPC-MIP, we list just the computation times. For CPC-MIPUB, we list the computation times and the speedup. For both CPC-LSCP and CPC-LSCPe we list the number of iterations, the computation times, and the speedup.

- **CPC-MIPUB** uniformly runs faster than CPC-MIP, and has an average speedup in our examples over CPC-MIP of 3.05, ranging from 1.74 to 10.22.
- **CPC-LSCP** uniformly runs faster than both of the MIP formulations, with an average speedup in our examples over CPC-MIP of 48.82 ranging from 2.61 to 369.89.
- **CPC-LSCPe** with enumeration for $p = 2$ and 3 uniformly runs faster than without, with an average speedup in our examples over CPC-MIP of 81.00 ranging from 3.06 to 671.05.
- **LSCP Iterations:** CPC-LSCPe cuts down on an average of 45% of the iterations relative to CPC-LSCP, ranging from a 17% cut to a 59% cut. Even though it is using enumeration to solve for only two values of $p$, this represents a major proportion of iterations because typically the covering distance gap for consecutive small $p$ is larger than the distance gap for consecutive large $p$ (see Fig. 1).
- Overall, CPC-LSCPe is the fastest of our methods for solving the CPC problem when using real number distances.

It is obvious, but worth noting nonetheless, that computation times increase exponentially and become quite large as the problem sizes increase. This is due to two primary factors. Firstly, any time one solves an NP-Hard problem on larger data sets, the computation times by definition increase non-polynomially. Secondly, as the problem size increases for the CPC problem, there are more PC subproblem solutions that must be found. Using any of the MIP approaches, a problem of size $n$ requires solving $n$ PC problems. Using the LSCP methods, as one increases the problem size, the set of unique distances in the distance matrix increase in addition to each iteration operating over a larger problem. This presents computational difficulties and limits the size of problems currently able to be solved in a reasonable amount of time. Future research directions listed in the conclusion may help to enable solution on larger problem sizes.

## 4.3 Experimental Evaluation Comparing Real Number to Integer Distances

We compared solution times on many of the TSP Library data sets for both real number distances and integer distances. Computing with the MIP algorithms, there was no real difference in computation times between problem variants, in-fact the integer instances were very slightly slower because one has to calculate the Euclidean distance matrix rather than the squared-Euclidean distance matrix. The significant differences came when solving the two problem variants with the LSCP algorithms. Rounding distances to integers significantly reduces the number of unique distances in the distance matrix, and thus the integer problems CPCi-LSCP and CPCi-LSCPe could complete far more quickly than their respective real number variants. Timed results are listed in Table 2, speedups compare an integer method to its respective real number method.

**Table 2** Computational experiments comparing real number distance CPC algorithms with integer distance CPCi algorithms.

| Name | Size | CPC-LSCP | CPCi-LCSCP | | CPC-LSCPe | CPCi-LSCPe | |
|---|---|---|---|---|---|---|---|
| | | Time (s) | Time (s) | Speedup | Time (s) | Time (s) | Speedup |
| rat195.tsp | 195 | 16.80 | 0.86 | 19.53 | 11.71 | 1.24 | 9.46 |
| d198.tsp | 198 | 18.89 | 5.91 | 3.20 | 7.61 | 2.33 | 3.26 |
| gr202.tsp | 202 | 127.09 | 0.30 | 427.58 | 108.31 | 0.79 | 136.39 |
| tsp225.tsp | 225 | 48.68 | 1.30 | 37.36 | 25.96 | 1.61 | 16.10 |
| gil262.tsp | 262 | 38.27 | 1.30 | 29.39 | 15.62 | 1.79 | 8.74 |
| a280.tsp | 280 | 8.27 | 1.44 | 5.73 | 6.36 | 2.74 | 2.32 |
| pr299.tsp | 299 | 62.38 | 21.88 | 2.85 | 34.85 | 15.55 | 2.24 |
| lin318.tsp | 318 | 113.64 | 20.07 | 5.66 | 50.40 | 13.00 | 3.88 |
| gr431.tsp | 431 | 1491.93 | 2.70 | 553.58 | 973.80 | 4.07 | 239.28 |
| pr439.tsp | 439 | 234.47 | 96.40 | 2.43 | 103.72 | 55.36 | 1.87 |
| u574.tsp | 574 | 3694.59 | 57.10 | 64.70 | 2453.82 | 47.03 | 52.18 |
| p654.tsp | 654 | 106.02 | 40.77 | 2.60 | 58.44 | 30.96 | 1.89 |
| fl1400.tsp | 1400 | 48110.86 | 233.42 | 206.11 | 10673.91 | 197.11 | 54.15 |

Integer distance problems always solved more quickly than their real number variants because they have a smaller set of unique distances to iterate over. There were a wide range of speedups, but CPCi-LSCP solved up to 553 times faster than the real number CPC-LSCP. On the largest problem we tested, with 1,400 points, the real number variant took 48,110 seconds to solve for all 1,400 $p$-values, while the integer variant took 233 seconds to solve. Adding enumeration to the integer problem, for small problems the number or unique distances in the distance matrix were often reduced enough that it was not beneficial to incorporate enumeration. For larger problems though it still remained advantageous. In the 1,400 point integer distance problem, CPCi-LSCPe reduced computation to 197 seconds as compared to 233 seconds without enumeration.

Rounding distances up and solving CPCi problem provides an integral feasible upper bound to the CPC real number distance problem. Subtracting one from those solutions provides an integral lower bound. We think future research could use this property to significantly narrow the search window and speed up computations for real number distance CPC solutions. This is only valid by rounding up (or down), rather than rounding to the nearest integer.

## 5 Conclusions

We introduced the complete *p*-center problem, and how it could be useful to planners and decision makers in determining the optimal trade-off between the number of facilities and the maximum distance from any demand to such a facility. We developed a basic formulation for the problem iterating the classical *p*-center MIP formulation over all *m* potential *p*-values, called CPC-MIP. We then showed a modification called CPC-MIPUB that trims many variables and constraints of the CPC-MIP formulation using an upper bound found from previous solutions while still maintaining optimality, resulting in an average speedup of 3.05.

Next, we developed a new algorithm called CPC-LSCP that breaks the CPC problem into many more LSCP sub-problems, solving for coverage with the unique distances in the distance matrix. While this significantly increases the number of iterations, each LSCP model has much fewer variables and constraints than the PC-MIP formulations and its integer friendly properties make it solve extremely quickly, resulting in an average speedup over the CPC-MIP problem of 48.82. Finally, we proposed the CPC-LSCPe algorithm by adding enumeration for $p = 2$ and 3 to the CPC-LSCP method, allowing it to solve with much fewer LSCP iterations. Enumeration was accomplished with highly optimized parallel compiled code, resulting in an average overall speedup for the CPC problem of 81.00 compared to CPC-MIP. While the exact performance numbers will likely vary by machine, the overall trend shows strongly that CPC-LSCPe is the fastest of these methods for solving the CPC problem for real number distances.

We also compared real number problem variants with integral variants, and found that integral distance CPCi problems solved much more quickly using the CPCi-LSCPe method on large problems due to having much fewer unique distances in their distance matrices. Every application must determine the precision unto which it must be solved, but the integer distance variant is undoubtedly easier to compute.

The complete *p*-center problem presents numerous topics for future research, including: 1) Solving the partial cover *p*-center problem to add isolines of partial coverage to the *p* vs. $z_p$ figure (see Fig. 1) to provide additional information to a designer or planner. 2) Speeding up the PC-LSCP method by skipping unique distances in the distance matrix, then backing up when there is a change to hone in on the exact coverage distance. Skipping every other distance could cut iterations by almost ½, higher skipping could potentially find larger speedups, binary search could result in logarithmic speedups. 3) Solving the integer distance CPCi problem with all distances rounded up provides an integral upper and lower bound for the real number version of the problem. The CPCi problem can be solved much more quickly than the real number distance variant, and these bounds would significantly reduce the search space for LSCP search of real number solutions. 4) We have demonstrated that methods can be combined for different *p*-values in order to optimize performance across all *p*. CPC-LSCPe uses the location set covering problem to quickly find solutions for high *p*-values where distance ties allow multiple solutions to be found with a single iteration, then uses enumeration for $p = 2$ and 3, and matrix math for $p = 1$ to quickly solve for low *p*-values. For solving the CPCi integer problem, a method could use CPCi-LSCP for high *p*-values then switch to CIK2018 for smaller *p* such as 30 and below. A faster integer method could also assist in solving the real number variant.

The complete *p*-center problem opens up several new and important directions in facility location research. It provides a planning tool that is rich with information for decisionmakers to weigh alternatives for facility location, and motivates the development of new algorithms that will be more efficient across the entire spectrum of *p*-values for the *p*-center problem. These research extensions will add value to decision sciences, OR algorithm development, and the facility location literature.

## Acknowledgments

We thank the Editor-in-Chief, guest editors, and two anonymous referees whose extensive comments and suggestions greatly helped to improve this manuscript. We also thank Greg Glockner and Gurobi for helping to optimize our code.

## References

Balinski ML (1965) Integer programming: methods, uses, computations Management science 12:253-313
Calik H, Labbé M, Yaman H (2015) p-Center problems. In: Location Science. Springer, pp 79-92
Calik H, Tansel BC (2013) Double bound method for solving the p-center location problem Computers & operations research 40:2991-2999
Chen D, Chen R (2009) New relaxation-based algorithms for the optimal solution of the continuous and discrete p-center problems Computers & Operations Research 36:1646-1655
Church R, Medrano F (2018) AM-46-Location-allocation modeling The Geographic Information Science & Technology Body of Knowledge doi:10.22224/gistbok/2018.3.4
Church R, ReVelle CR (1974) The maximal covering location problem Papers in Regional Science 32:101-118
Church R, Weaver J (1986) Theoretical links between median and coverage location problems Annals of Operations Research 6:1-19
Church RL (2008) BEAMR: An exact and approximate model for the p-median problem Computers & Operations Research 35:417-426
Church RL (2016) Tobler's Law and Spatial Optimization Why Bakersfield? International Regional Science Review
Church RL, Gerrard RA (2003) The multi-level location set covering model Geographical Analysis 35:277-289
Church RL, Roberts KL (1983) Generalized coverage models and public facility location Papers of the Regional Science Association 53:117-135 doi:10.1007/bf01939922
Contardo C, Iori M, Kramer R (2019) A scalable exact algorithm for the vertex p-center problem Computers & Operations Research 103:211-220
Daskin MS (2013) Center Problems. In: Network and Discrete Location: Models, Algorithms, and Applications. 2nd edn. John Wiley & Sons, Inc., pp 193-234. doi:10.1002/9781118537015.ch05
Daskin MS, Maass KL (2015) The p-median problem. In: Location science. Springer, pp 21-45
Daskin MS, Owen SH (1999) Two New Location Covering Problems: The Partial P-Center Problem and the Partial Set Covering Problem Geographical Analysis 31:217-235
Drezner Z (1984) The planar two-center and two-median problems Transportation Science 18:351-361
Dzator M, Dzator J (2013) An effective heuristic for the P-median problem with application to ambulance location Opsearch 50:60-74
Elloumi S, Labbé M, Pochet Y (2004) A new formulation and resolution method for the p-center problem INFORMS Journal on Computing 16:84-94
Fo A, da Silva Mota I (2012) Optimization models in the location of healthcare facilities: a real case in Brazil Journal of Applied Operational Research 4:37-50
Gadegaard SL, Klose A, Nielsen LR (2018) An improved cut-and-solve algorithm for the single-source capacitated facility location problem EURO Journal on Computational Optimization 6:1-27
García S, Labbé M, Marín A (2011) Solving large p-median problems with a radius formulation INFORMS Journal on Computing 23:546-556
Hakimi SL (1964) Optimum locations of switching centers and the absolute centers and medians of a graph Operations research 12:450-459

Kariv O, Hakimi SL (1979) An algorithmic approach to network location problems. I: The p-centers SIAM Journal on Applied Mathematics 37:513-538

Marianov V, ReVelle C (1995) Siting emergency services Facility Location: a survey of applications and methods 1:199-223

Medrano FA (2019) p-center. https://github.com/antoniomedrano/p-center. Accessed September 25, 2019

Minieka E (1970) The m-center problem Siam Review 12:138-139

Mladenović N, Labbé M, Hansen P (2003) Solving the p-center problem with tabu search and variable neighborhood search Networks 42:48-64

Niblett MR, Church RL (2015) The disruptive anti-covering location problem European Journal of Operational Research 247:764-773

Reinelt G (1991) TSPLIB—A traveling salesman problem library ORSA journal on computing 3:376-384

ReVelle C (1993) Facility siting and integer-friendly programming European Journal of Operational Research 65:147-158

Rosing KE, ReVelle C, Rosing-Vogelaar H (1979) The p-median and its linear programming relaxation: An approach to large problems Journal of the Operational Research Society 30:815-823

Salazar-Aguilar MA, Ríos-Mercado RZ, González-Velarde JL (2013) GRASP strategies for a bi-objective commercial territory design problem Journal of Heuristics 19:179-200

Schilling DA (1982) Strategic facility planning: The analysis of options Decision Sciences 13:1-14

Shier D (1977) A min-max theorem for p-center problems on a tree Transportation Science 11:243

Snyder LV, Daskin MS (2006) Stochastic p-robust location problems IIE Transactions 38:971-985

Solanki R (1991) Generating the noninferior set in mixed integer biobjective linear programs: an application to a location problem Computers & Operations Research 18:1-15

Suzuki A, Drezner Z (1996) The p-center location problem in an area Location science 4:69-82

Swain RW (1971) A decomposition algorithm for a class of facility location problems.

Teitz MB, Bart P (1968) Heuristic methods for estimating the generalized vertex median of a weighted graph Operations Research:955-961

Tutunchi GK, Fathi Y (2019) Effective methods for solving the Bi-criteria p-Center and p-Dispersion problem Computers & Operations Research 101:43-54

Wei H (2008) Solving Continuous Space Location Problems. Dissertation, The Ohio State University